\documentclass[11pt,a4paper]{amsart}
\usepackage{amscd,amssymb,amsmath,latexsym,enumerate}

\newcommand{\cha}{{\operatorname{char}}}

\newcommand{\bfk}{{\boldsymbol{k}}}

\newcommand{\bfzero}{{\boldsymbol{0}}}

\theoremstyle{plain}

\newtheorem{theorem}{Theorem}[section]

\newtheorem{proposition}[theorem]{Proposition}
\newtheorem{corollary}[theorem]{Corollary}

\newtheorem{definition}[theorem]{Definition}

\newtheorem{remark}[theorem]{Remark}

\def\N{{\mathbb N}}

\def\Z{{\mathbb Z}}

\begin{document}

\title[Multivariate symmetric Hermite Interpolant]{A note on the multivariate \\symmetric Hermite Interpolant}

\author[T. Krick]{Teresa Krick}
\address{Departamento de Matem\'atica, FCEN, Universidad de Buenos Aires and IMAS, UBA-CONICET,  Argentina.}
\email{krick@dm.uba.ar}
\urladdr{http://mate.dm.uba.ar/\~\,krick}

\author[A. Szanto]{Agnes Szanto}
\address{Department of Mathematics, North Carolina State
University, Raleigh,  USA.}
\urladdr{https://aszanto.math.ncsu.edu/}

\begin{abstract} In this note we explicit the  notion of Hermite interpolant of a multivariate symmetric polynomial, generalizing the notion of Lagrange interpolant to the case when there are roots coalescence, an extension of the results on the symmetric Hermite interpolation basis by M.-F. Roy and A. Szpirglas included in \cite{RoSz2020}.
\end{abstract}


\maketitle

\section{Introduction}

{\em We started this project with Agnes  in 2017 but never found the occasion to develop it in depth.  In 2022 Agnes passed away. I finally found the time and emotional calm to give an end to this project. This is thanks  to Agnes, my very special friend and greatest collaborator, and for her (T.K).}

\bigskip
We are all familiar with   Lagrange interpolation which says that for  $d\in \N$, given any set  $A$ with $d$ elements in some field $K$, i.e. $|A|=d$, the set of $d$  univariate polynomials
$$\Big\{\frac{\prod_{a'\ne a} (x-a')}{\prod_{a'\ne a} (a-a')}\ : \quad a\in A\Big\}$$
is a basis, the Lagrange basis, of the vector space of all polynomials in $K[x]$  of degree $< d$. \\
In particular, for any choice of $d$ constants $c_{a}\in K, \, a\in A$, the  unique polynomial $r\in K[x]$ of degree $<d$ which satisfies the $d$ conditions $r(a)=c_a$ for all $a\in A$  is the polynomial
$$r= \sum_{a\in A} c_a \frac{\prod_{a'\ne a} (x-a')}{\prod_{a'\ne a} (a-a')}.$$ For any polynomial $h\in K[x]$ of arbitrary degree, the  unique polynomial $r_A(h)\in K[x]$ of degree $< d$ which satisfies the $d$ conditions  $r_A(h)(a)=h(a)$ for all $a\in A$, which can be expressed in this basis as
$$r_A(h)=\sum_{a\in A} h(a) \frac{\prod_{a'\ne a} (x-a')}{\prod_{a'\ne a} (a-a')},$$
is called the {\em  Lagrange interpolant}   of $h$ with respect to the set $A$.
\\
It is clear that if we associate to the set $A$ the polynomial $f=\prod_{a\in A} (x-a) \in K[x]$ whose roots are the elements of $A$, then $r_A(h)=r_f(h)$, where $r_f(h)$ denotes the remainder of $h$ when dividing it by $f$. 

\smallskip
Less known is the fact that Lagrange interpolation extends naturally to  {\em symmetric polynomials in  $n$ variables} $X=\{x_1,\dots,x_n\}$ of degree in each  variable  $ \le d-n$. We were thrilled when in 2014  we  (re)discovered this result  and showed its impact for our research related to subresultants and Sylvester's double sums in \cite{KrSzVa2017} and \cite{DKSV2019}. But then   we found that it had already appeared (as novel) in 1996 in    \cite{ChLo1996}, and before that  in  \cite{Bor1860} and even  in the appendix of  Jacobi's thesis in 1825 \cite{Jac1825}!
\\
Before stating the symmetric Lagrange interpolation, we set the following notations.
Given any finite sets $Y,Z$,  we denote
$$R(Y,Z):=\prod_{y\in Y, z\in Z} (y-z)$$
with $R(Y,Z)=1$ if $Y=\emptyset$ or $Z=\emptyset$ (the notation $R$, following A. Lascoux, is chosen because  it coincides with the resultant $R(f,g)$ of the univariate polynomials $f$ and $g$ with roots in $Y$ and $Z$ respectively). If $Y=\{y\}$ we simply write $R(y,Z)$.
Also, given $0\le n\le |Y|$, we denote by  $Y' {\subset_n} Y$ the fact that  $Y'$ is a subset of $Y$ with exactly  $n$ elements, and we note that there are exactly  $\binom{|Y|}{n}$ such subsets $Y'$.\\
Finally, given the set of  $n$ variables $X=\{x_1,\dots,x_n\}$, we denote by $S[X]\subset K[X]$ the $K$-algebra of {symmetric} polynomials in $K[X]$, and  for $d\in \N$ with $d\ge n$,  we denote by  $S_{d-n}[X]\subset S[X]$  the $K$-vector space of symmetric polynomials in $S[X]$ of degree in each variable  $\le d-n$.
\\
We are now ready to  state the theorem, simple proofs of it can be found for instance in \cite[Th.2.1]{ChLo1996}  and \cite[Prop.2.3]{KrSzVa2017}.

\begin{theorem}\label{thm:Lagrange} {\em (Symmetric Lagrange interpolation.)}\\
Set  $n\le d\in \N$. Let $A\subset  K$ with  $|A|=d$ and $X=\{x_1,\dots,x_n\}$. Then the set of $\binom{d}{n}$ symmetric polynomials:
$$\Big\{ \frac{R(X,A\setminus A')}{R(A',A\setminus A')}=\frac{\prod_{{x_i\in X,a'\in A\setminus A'}} (x_i-a')}{\prod_{{a\in A',a'\in A\setminus A'}} (a-a')}\ :  \quad A'{\subset_n}A\Big\}\ \subset \  S_{d-n}[X]$$
is  a basis, the {\em Lagrange basis}, of the vector space $S_{d-n}[X]$ of all symmetric polynomials in  $n$ variables of degree in each variable $\le d-n$.

In particular, for any choice of $\binom{d}{n}$ values $c_{A'}\in K, \, A'\subset_n A$, the  unique symmetric polynomial $r\in S_{d-n}[X]$ of degree $\le d-n$ in each variable which satisfies the $\binom{d}{n}$ conditions $r(A')=c_{A'}$ for all $A'\subset_n A$  is the polynomial
$$r= \sum_{A'{\subset_n}A} c_{A'} \frac{R(X,A\backslash A')}{{{R}(A',
A\backslash A')}}.$$
(Here evaluating a symmetric polynomial into $A'$ means evaluating its variables into the elements of $A'$ in any order, which is well-defined since the polynomial is symmetric.)\end{theorem}
This allows us to give the following definition.
\begin{definition} \label{def:interpolant}{\em (Symmetric Lagrange interpolant.)}\\
Given a set $A\subset K$ with $|A|=d\ge n$ and given a symmetric polynomial $h\in S[X]$, the {\em symmetric Lagrange interpolant} in  $S_{d-n}[X]$ of $h$ with respect to the set $A$ is the unique symmetric polynomial $r_A(h)\in S_{d-n}[X]$ which  satisfyes the $\binom{d}{n}$ conditions $r_A(h)(A')=h(A')$  for all $A' \subset_n A$. It is given by the expression
\begin{eqnarray*}\label{Lagformula}r_A(h) = \sum_{A'{\subset_n}A} h(A') \frac{R(X,A\backslash A')}{{{R}(A',
A\backslash A')}}.
\end{eqnarray*}
\end{definition}
Note that in the multivariate case,  defining $r_A(h)$ as a remainder of a polynomial  division in $K[X]$ is not so clear and this is one of the problems we address in this note.

\smallskip
Now, the (univariate) Lagrange interpolation classically extends  to  Hermite interpolation, which is the case when the nodes in $A$ can be repeated, i.e. $A$ is a multiset
\begin{equation}\label{eq:multiset} A=
\big\{\underbrace{a_1,\dots,a_1}_{d_1},\dots,\underbrace{a_m,\dots,a_m}_{d_m}\big\} \ \subset \, K\end{equation} with $|A|=d:=d_1+\dots+ d_m$: There is a unique    (univariate) polynomial $r\in K[x]$ of degree $<d$  which satisfies the $d$ conditions
\begin{equation}\label{eq:Hermite}\left\{\begin{array}{l}r(a_1)=c_{1,0},\ r'(a_1)=c_{1,1},\ \dots\ ,\ r^{(d_1-1)}(a_1)=c_{1,d_1-1},   \\
\ \ \vdots \\
r(a_m)=c_{m,0},\ r'(a_m)=c_{m,1},\ \dots \ ,\ r^{(d_m-1)}(a_m)=c_{m,d_m-1}, \end{array}\right.
\end{equation}
for given $c_{1,0},\dots,c_{1,d_1-1},\dots,c_{m,1},\dots,c_{m,d_m-1} \in K$.
\\
For Hermite interpolation there is no simple  description of the natural Hermite basis $\big\{r_{i,k_i}\in K[x], \ 1\le i\le m, 0\le k_i<d_i\big\}$,
which satisfies that
$$r=\sum_{i,k_i}c_{i,k_i} r_{i,k_i},$$  as in the Lagrange case, except maybe for  the case of the Taylor interpolation where there is a single node $a\in A$ and the basis is given by $$\Big\{\frac{(x-a)^k}{k!}\,:\ 0\le k<d\Big\}.$$ However, there is  a simple matrix description of this basis that comes from the generalized Vandermonde (or confluent) matrix which describes the Hermite linear  problem~\eqref{eq:Hermite} (see Section \ref{sec:RS}  below for more details). \\
We also remark that the Hermite interpolant $r_A(h)$ of an arbitrary polynomial $h\in K[x]$, which is the unique polynomial of degree $<d$ that satisfies the $d$ conditions in \eqref{eq:Hermite} for $c_{k,i_k}=h^{(i_k)}(a_k)$, is again easily recovered as
$r_A(h)=r_f(h)$, the remainder of $h$ when dividing it by the polynomial
\begin{equation}\label{eq:f}f=\prod_{1\le i\le m}(x-a_i)^{d_i} \quad \in \ K[x]\end{equation} associated to the multiset $A$.

\smallskip

In \cite{RoSz2020},  M.-F.Roy and A. Szpirglas present, as a step for their description  of Sylvester's double sums for polynomials with multiple roots, an extension of  the symmetric Lagrange basis of Theorem~\ref{thm:Lagrange} to the case when  $A$ is a multiset. More precisely,  they describe a basis $\mathcal{B}$ for the multivariate symmetric Hermite interpolation problem, which agrees with the Hermite basis in the univariate case and with the multivariate symmetric Lagrange basis when $A$ is a set.   Their construction, as explained in Section~\ref{sec:RS} below,  is  a matrix construction (as in the univariate Hermite case), which  describes the polynomials in this basis $\mathcal{B}$ as determinants of matrices composed by certain columns  of generalized Vandermonde matrices associated to the variables in $X$ and the multiset $A$. They also describe the coordinates of a given symmetric polynomial $r\in S_{d-n}[X]$ in $\mathcal{B}$, by means of ordered evaluations in $A$ of appropriate derivatives   of $v_n\, r$ where $v_n=\prod_{j>i}(x_j-x_i)$ is the Vandermonde determinant of the (ordered) variables in $X$.  

\smallskip
In this note we go one step further and first provide  a natural extension  of the notion of Lagrange interpolant  $r_A(h)\in S_{d-n}[X]$ of any symmetric polynomial $h\in S[X]$ introduced in Definition~\ref{def:interpolant} to the case of  a multiset, which we call the {\em Hermite interpolant} of $h$. We also show that  the same conditions on the derivatives of $v_n\, r$ that describe the coordinates in the basis  $\mathcal{B}$  of  any $r\in S_{d-n}[X]$ hold for the derivatives of $v_n\, h$ to describe the coordinates of $r_A(h)$ in $\mathcal{B}$ of any symmetric polynomial $h\in S[X]$.

\smallskip
Before stating our main result,  we need to introduce a special family of polynomials. First, let the polynomial $f$ in \eqref{eq:f}, whose roots are the elements in the multiset $A$ with the corresponding multiplicities, be written as
$$f=x^d+f_{d-1}x^{d-1}+\cdots + f_1x+f_0 \ \in \, K[x],$$
where $f_{d-k}=(-1)^k s_k(A)$, the signed elementary symmetric polynomial   $s_k$ of degree $k$ in $d$ variables evaluated into the elements of the multiset $A$.  \\
Next, inspired by \cite[Prop.2.1\& Thm.2.2]{HeNaRo2004}, we introduce 
the following family $G=\{g_1\dots,g_n\}\subset K[X]$ of $n$ polynomials in $n$ variables defined from the coefficients $f_{d-k}$, $1\le k\le d$,  of $f$:
\begin{equation}\label{GB}
\begin{array}{cc}G:\left\{\begin{array}{l}
g_1(x_1):=f(x_1) = h^{(1)}_d+ f_{d-1}h_{d-1}^{(1)}+\cdots+f_1h_{1}^{(1)}+ f_0 ,\\[1mm]
g_2(x_1, x_2):=h_{d-1}^{(2)}+ f_{d-1}h_{d-2}^{(2)}+\cdots+f_2h_{1}^{(2)} + f_1,  \\[1mm]
 \quad\vdots \\
g_{n}(x_1, \ldots, x_{n}):=h_{d-n+1}^{(n)}+ f_{d-1}h_{d-n}^{(n)}+\cdots+f_{n-1},\end{array}\right.\end{array}
\end{equation}
where $h_j^{(i)}$ denotes the  {\em complete homogeneous symmetric polynomial}  of degree $j$ in $x_1, \ldots ,x_i$, which is the sum of all monomials of total
degree $j$ in $x_1,\dots,x_i$:
\begin{eqnarray}\label{complete}
h_{j}^{(i)}:=\sum_{|\bfk|=j}x_1^{k_1}\cdots x_i^{k_i}, \ 1\le i\le n, \  j\ge 0.
\end{eqnarray}
Given $1\le i\le n$, the  polynomial $g_i$ is symmetric by construction, but moreover,
we will see in Section~\ref{subsec:hermite} that the ideal $(G)\subset K[X]$ generated by the polynomials $g_1,\dots,g_n$ in $G$ is a symmetric ideal (i.e., if a polynomial $h(X)\in (G)$, then $h(X^\sigma)\in (G)$ for any permutation $\sigma$ of the variables  in $X$). Therefore since $f(x_1)=g_1\in (G)$, $f(x_i) \in (G)$ for $1\le i\le n$. More precisely, we show  in the proof of Proposition~\ref{prop:Vand} that  the following relationship holds:
\begin{equation}\label{eq:fg}
\left\{\begin{array}{l}
f(x_1) =g_1\\
f(x_2)= (x_2-x_1)\,g_2+g_1\\
f(x_3)=(x_3-x_1)(x_3-x_2)\,g_3+ (x_3-x_1)\,g_2+g_1\\
\quad \vdots\\
f(x_{n}) = \prod_{i=1}^{n-1}(x_{n}-x_i)\, g_{n} + \cdots + (x_{n}-x_1)\,g_{2}+g_{1}.\end{array}\right.
\end{equation}
Note that each  of the polynomials $g_i$ belongs to $K[x_1,\dots,x_i]$ and is monic of degree $d-i+1$ in the variable $x_i$.  Given a  polynomial $h\in K[X]$, it then makes   sense to
consider $r_G(h)\in K[X]$, the remainder of dividing $h$ by the family $G$, which  corresponds  to performing successive  Euclidean divisions, first of $h$ by  the polynomial  $g_n$, monic in $x_n$,  in $\big(K[x_1,\dots,x_{n-1}] \big)[x_n]$, second  divide its coefficients when viewed as a polynomial in $x_n$ by $g_{n-1}$ in $\big(K[x_1,\dots,x_{n-2}] \big)[x_{n-1}]$, and so on. This is formalized by  taking the normal form of $h$ with respect to the Gr\"obner basis $G$ for the pure lexicographic order $\prec$ with $x_1\prec x_2\prec \dots \prec x_n$,  as explained in Proposition~\ref{grobner} below.

\smallskip
We note that in the case of a single variable $x$, $G=\{f\}$, and therefore $r_G(h)=r_f(h)$ is the Hermite interpolant of the polynomial $f$.
Our main result generalizes this fact to the multivariate symmetric Hermite interpolant of a multivariate symmetric polynomial.
\begin{theorem}\label{thm:main} {\em (Symmetric Hermite interpolant.)}\\
Let   $A\subset K$ be a multiset as in \eqref{eq:multiset} with $|A|=d\ge n$ and let $h\in S[X]$ be a symmetric polynomial.
 Then,    $$r_G(h)\in S_{d-n}[X].$$
  Moreover,  $r_G(h)$  is a polynomial expression in the elements of $A$ which coincides with the Hermite interpolant  in the univariate case and with the symmetric Lagrange interpolant  introduced in Definition \ref{def:interpolant} in the  case of a set $A$.\\ Thus, the symmetric Hermite interpolant  $r_A(h)$ is equal to $r_G(h)$ for a multiset $A$ and $h\in S[X]$.
\end{theorem}
Our second result is presented in Corollary \ref{cor:Hermiteinterpolant} below and describes the coordinates of the Hermite interpolant $r_G(h)$ of any symmetric polynomial $h\in S[X]$ in terms of the symmetric Hermite basis $\mathcal{B}$  introduced in \cite{RoSz2020}.

\smallskip

The structure of this note is the following: Section~\ref{sec:main}  proves Theorem \ref{thm:main}. In Subsection~\ref{subsec:lagrange} we first consider the Lagrange case where all the points in $A$ are different and show in Proposition~\ref{grobner} that in that case $r_A(h)=r_G(h)$. We then show in Subsection~\ref{subsec:hermite} that in the  Hermite case, when they are points coalescence, $r_G(h)$ is a polynomial in the elements of $A$ that belongs to $ S_{d-n}[X]$ for any $h\in S[X]$, allowing to finish the proof of Theorem~\ref{thm:main}. Section~\ref{sec:RS} makes the connection between the Hermite interpolant and Roy-Szpirglas' Hermite basis: We begin in Subsection~\ref{sec:RS2} by recalling the construction in \cite{RoSz2020}, and in Subsection~\ref{sec:coord} we present our second main result, Corollary~\ref{cor:Hermiteinterpolant}. We end this note by Section~\ref{sec:taylor}, which illustrates the results for symmetric Taylor interpolation,  i.e. the  case of a multiset $A$ with only one repeated node.

\bigskip
\noindent {\bf Acknowledgment.} T.K. is grateful to Nicol\'as Allo G\'omez,  Mat\'\i as Bender, Marie-Fran\c coise Roy, Aviva Szpirglas , and to the anonymous referees, for their comments, which greatly helped improve the presentation of this note.

\section{Proof of Theorem~\ref{thm:main}}\label{sec:main}

We recall that  $K$  is a field, $1\le n\le d$ in $\N$, $X=\{x_1,\dots,x_n\}$, $S[X]\subset K[X]$ is the algebra of symmetric polynomials in the $n$ variables in $X$ and $S_{d-n}[X]$ is the vector space of polynomials in $S[X]$ with degree in each variable $\le d-n$.

\subsection{The symmetric Lagrange interpolant}\label{subsec:lagrange}{\ }

 \smallskip
In this section $K$ can be of any characteristic,  $A\subset K$ is a set of cardinality $|A|=d$ and 
$$f=\prod_{a\in A} (x-a) \ = \ x^d+f_{d-1}x^{d-1}+\cdots + f_1x+f_0 \quad \in \ K[x] $$
is as previously the (squarefree) univariate monic polynomial whose   roots are the (distinct in this case) elements in $A$. 
We introduce the  following symmetric set in $A^n$ of $n$-tuples with distinct elements of $A$
$$
V_n(A):= \left\{ \left(a_{1}, \ldots, a_{{n}}\right) \ : \ a_1,\dots,a_n\in A,\  a_i\ne a_j \ \mbox{for} \ 1\le i\ne j\le n\right\},
$$ which satisfies  $\big| V_n(A)\big|=d!/(d-n)!$ and
we  set
$$
I:=I(V_{n}(A))\subset K[X]
$$ for its  vanishing ideal in $K[X]$, which is a zero-dimensional radical  ideal in $K[X]$, which is moreover a symmetric ideal,   since it is the vanishing ideal of a finite symmetric set.

\smallskip
Next proposition relates the symmetric Hermite interpolant $r_A(h)$ introduced in Definition~\ref{def:interpolant} of  a symmetric polynomial $h$ with the remainder $r_G(h)$ of $h$ with respect to the family $G$ introduced in \eqref{GB}.
Its first claim  appears for the $n=d$ case in \cite[Thm.2.2]{HeNaRo2004}.

\begin{proposition}\label{grobner} {\ }
\begin{enumerate}
\item
The set
$$
G :=\left\{ g_1(x_1), g_2(x_1, x_2), \dots,
g_{n}(x_1, \ldots, x_{n})\right\}
$$ defined  in  \eqref{GB}
is a Gr\"obner basis for $I$ with respect to  the lexicographic order $\prec$ of $K[X]$ with $x_1<x_2<\cdots<x_{n}$, and the normal set of
$K[X]/I$  with respect to $\prec$  is spanned by the basis
$$B:=\{ x_1^{k_1}\cdots x_{n}^{k_{n}}:0\le k_1\le d-1,\dots,0\le k_{n}\le d-n\}.$$
\item
Given a {\em symmetric} polynomial $h\in K[X]$, then
$$r_A(h)=r_G(h),$$
where $r_G(h)$ denotes the remainder (i.e. normal form) of $h$ modulo the Gr\"obner basis $G$ with respect to the order $\prec$.
\end{enumerate}
\end{proposition}

\begin{proof}
$\it{(1)}$
It is clear that $G$ is a Gr\"obner basis of the ideal $( G)\subset K[X]$ it generates w.r.t. $\prec$, since the leading term of each $g_i$ equals $ x_i^{d-i+1}$, $1\le i\le n$, and they are  pairwise relatively prime.\\
Note that the  factor ring $K[X]/( G)$ is then spanned over $K$ by $$B=\{ x_1^{k_1}\cdots x_{n}^{k_{n}}:0\le k_1\le d-1,\dots,0\le k_{n}\le d-n\}.$$ Therefore
$$\dim_K\big(K[X]/( G))=  \frac{d!}{(d-n)!}=\#(V_{n}(A))=\dim_K\big(K[X]/I\big),$$ where the last equality holds because $I$ is a radical ideal.
  Thus it is enough to prove that  $G\subset I$, since in that case we have the following isomorphisms of finite-dimensional $K$-vector spaces
   $$ K[X]/( G) \simeq_K K[X]/I \simeq_K \big(K[X]/( G)\big)/\big(I/( G)\big)$$ which implies that $I=( G)$.  \\This is equivalent to show that $g_i$  vanishes on  $V_{n}(A)$, for all $i$.
  We show it by recurrence on $i$ for $i=1,\dots,n$:\\
It is clearly true for $g_1=f(x_1)$ since $f(a)=0$, $\forall a\in A$.\\
We now show that the case $i-1$ implies the case $i$: Using that for $i>1$,
$$
(x_1-x_i)h^{(i)}_{j}(x_1, \ldots, x_i) = h^{(i-1)}_{j+1}(x_1, \ldots, x_{i-1}) - h^{(i-1)}_{j+1}(x_i, x_2,\ldots, x_{i-1}),
$$
since the  monomials of degree $j+1$ that contain both $x_1$ and $x_i$ are canceled on the left hand side, and the monomials which do not contain any of $x_1$ and $x_i$ are canceled on the right hand side,  we derive that
\begin{equation}\label{bi}
g_i(x_1, \ldots, x_i)=  \frac{g_{i-1}(x_1, \ldots, x_{i-1})-g_{i-1}(x_2, \ldots, x_{i})}{x_1-x_i}.
\end{equation}
and therefore $g_i$ vanishes on  $V_{n}(A)$ since $g_{i-1}$ does and $a_1\ne a_i$ in $V_{n}(A)$.

\smallskip
\noindent
$\it{(2)}$
Let $h\in S[X]$ be a symmetric polynomial. By Theorem~\ref{thm:Lagrange},  $r_A(h)$ is the unique polynomial in $S_{d-n}[X]$ which satisfies
$r_A(h)(a_1,\dots,a_n) = h(a_1,\dots,a_n)$ for all $\binom{d}{n}$ choices of $\{a_1,\dots,a_n\} \subset A$.
 Since $h-r_A(h)$ is symmetric and hence vanishes on all  $V_{n}(A)$, we have $h-r_A(h)\in I$. Therefore $h=q+r_A(h)$ for some $q\in I$, where $r_A(h)\in S_{d-n}[X]$ is a valid remainder modulo the Gr\"obner basis $G$  since all its monomials have degree  in each variable $\le d-n$, and thus belong to  ${\rm span}(B)$. This shows that $r_A(h)= r_G(f)$.
\end{proof}

 \subsection{The symmetric Hermite interpolant}\label{subsec:hermite}{\ }

 \smallskip
Here we have a multiset  $A$ as introduced in \eqref{eq:multiset}, of cardinality $d$,  which defines the univariate polynomial
$$f=\prod_{1\le i\le m}(x-a_i)^{d_i} =x^d+f_{d-1}x^{d-1}+\cdots+f_0 \ \in \ K[x] .$$
As before we have
$$G :=\left\{ g_1(x_1), g_2(x_1, x_2), \dots,
g_{n}(x_1, \ldots, x_{n})\right\}\subset K[X]
$$  where  $g_1,\dots, g_{n}$ are defined in Identity~\eqref{GB}, which is a  Gr\"obner basis of the ideal
$( G) \subset K[X]
$
it generates with respect to  the lexicographic order $\prec$ of $K[X]$ with $x_1<\cdots <x_{n}$. (Note that here we cannot  define yet an ideal  $I$ as the vanishing set of $V_n(A)$ because we don't have the finite symmetric set $V_n(A)$ as in the previous section.)

 \begin{proposition} \label{inv}  Let $G$ be defined as in \eqref{GB}. Then, \begin{enumerate} \item
 $(G)\subset K[X]$ is a symmetric ideal, i.e. if a polynomial $g(X)\in (G)$, then $g(X^\sigma)\in (G)$ for any permutation $\sigma$ of the variables in $X$.
 \item For any {\em symmetric} polynomial $h\in S[X]$,
$$
r_G(h) \in S_{d-n}[X],
$$
i.e. it is a symmetric polynomial in $K[X]$ of degree in each variable  $\le d-n$.
\end{enumerate}
  \end{proposition}

 \begin{proof} {\it (1)} It is enough to prove that for all $1\le i\le n$, we have $$g_i(x_{j_1}, \ldots, x_{j_i})\in (G)$$ for any distinct $1\leq {j_1}, \ldots, j_i\leq n$. First note that $g_i(x_{j_1}, \ldots, x_{j_i})$ is symmetric in $x_{j_1}, \ldots, x_{j_i}$ for all distinct $1\leq {j_1}, \ldots, j_i\leq n$ by construction. In particular $g_n(x_{j_1},\dots,x_{j_n})=g_n(x_1,\dots,x_n)\in (G)$ and moreover, $g_i(x_{j_1}, \ldots, x_{j_i})=g_i(x_1,\dots,x_i)\in (G)$ for any distinct $1\leq {j_1}, \ldots, j_i\leq i$.
\\
 We will now prove that for  all $2\le i\le n$, $g_i(x_{j_1}, \ldots, x_{j_i})\in (G)$ for any distinct $1\leq {j_1}, \ldots, j_i\leq n$  implies that $g_{i-1}(x_{k_1}, \ldots, x_{k_{i-1}})\in (G)$   for any distinct $1\leq {k_1}, \ldots, k_{i-1}\leq n$. \\We will do it from knowing that $g_{i-1}(x_{k_1}, \ldots, x_{k_{i-1}})\in (G)$ for any distinct $1\le k_1,\dots, k_{i-1}\le i-1$ and by moving one variable at a time. \\
By Identity~\eqref{bi}  we have
 $$
g_i(x_{k_1}, \ldots,x_{k_i})=  \frac{g_{i-1}(x_{k_1}, \ldots, x_{k_{i-1}})-g_{i-1}(x_{k_2}, \ldots, x_{k_i})}{x_{k_1}-x_{k_i}}
$$
i.e. \begin{equation}\label{eq:idg}g_{i-1}(x_{k_1}, \ldots, x_{k_{i-1}})=(x_{k_1}-x_{k_i})\,g_i(x_{k_1}, \ldots,x_{k_i})+ g_{i-1}(x_{k_2}, \ldots, x_{k_{i}})\end{equation}
for any distinct  $1\le k_1,\dots,k_{i}\le n$.
Taking $\{k_2,\dots,k_i\}=\{1,\dots, i-1\}$ and $k_1\in\{ i,\dots,n\}$ we deduce from the fact that $g_i(x_{k_1}, \ldots,x_{k_i})\in (G)$ and $g_{i-1}(x_{k_2}, \ldots, x_{k_{i}})\in (G)$ that $g_{i-1}(x_{k_1}, \ldots, x_{k_{i-1}})\in (G)$. What we did here is to replace the variable $x_{k_i}\in\{x_1,\dots,x_{i-1}\}$ by a variable in $\{x_i,\dots, x_n\}$. We proceed replacing variables one at a time and conclude in the end that $g_i(x_{k_1}, \ldots,x_{k_i})\in (G)$ for any distinct  $1\le k_1,\dots,k_{i-1}\le n$.

\smallskip
\noindent {\it (2)} Let $h \in S[X]$ be a symmetric polynomial. Then $h=g+r_G(h)$ with $g\in (G)$ a symmetric polynomial, and therefore $r_G(h)$ is a symmetric polynomial. This together with $\deg_{x_n}( r_G(h))\le d-n$ implies that $\deg_{x_i}( r_G(h))\le d-n$ for $1\le i\le n$, i.e. $r_G(h)\in S_{d-n}[X]$.
\end{proof}
In view of the previous result, it makes sense to propose that for a  multiset $A\subset K$ as in \eqref{eq:multiset}, the Hermite interpolant of a symmetric polynomial $h\in S[X]$ is defined as  $r_G(h)\in S_{d-n}[X]$. We will justify this by showing that the polynomial $r_G(h)$ is a polynomial in the elements of the multiset $A$. To this aim,
we first observe that for any polynomial $h\in R[X]$, where $R\subset K$ is an integral domain, $r_G(h)$ is in fact a polynomial in the coefficients of $f$ rather than a rational function in them. Since  the leading term of each $g_i\in G$ equals the monic monomial $ x_i^{d-i+1}$, we never need to perform divisions by leading coefficients when performing the division of a polynomial $h$ by a $g_i$.
 This allows us to pass to generic polynomials: 

\begin{remark}Let $U=\{u_{0},\dots,u_{d-1}\}$ be a set of $d$ distinct indeterminates and $R$ be an integral domain. Let $$u=x^d+u_{d-1}x^{d-1}+\cdots + u_0 \ \in \ \Z[U][x]$$ be the generic monic polynomial of degree $d$ with indeterminate coefficients in $U$.  Define the corresponding generic polynomials  $$g_i(U,x_1,\dots,x_i)=h_{d-i+1}^{(i)}+ u_{d-1}h_{d-i}^{(i)}+ \dots + u_{i-1} \ \in \ \Z[U][X]\quad \mbox{for} \ 1\le i\le n,$$
where $h_j^{(i)} $ is the complete homogeneous polynomial introduced in \eqref{complete}.
Then, since the leading term of each $g_i$  as a polynomials in $X$   w.r.t. the lexicographic order $\prec$ equals $x_i^{d-i+1}$, which are all pairwise relatively prime, the  set $$G(U):= \{ g_1(U,x_1),g_2(U,x_1,x_2),\dots , g_n(U,x_1,\dots,x_n)\}  \, \subset \, \Z[U][X]$$  is a Gr\"obner basis of the ideal $( G(U)) \subset \Z(U)[X]$ w.r.t. the lexicographic order $\prec$ of the variables $X$.\\
Moreover, for any  polynomial $h\in R[X]$, the fact that the leading coefficient of each of these generic  $g_i(U,X)$'s as polynomials in $X$ equals 1 implies that   $r_{G(U)}(h)\in R[U][X]$, i.e. its coefficients are polynomials in $U$ rather than rational functions.
\end{remark}
We derive the following corollary in terms of a polynomial $u$ determined by  generic roots:
\begin{corollary} Let $Y=\{y_1,\dots,y_d\}$ be a set of $d$ distinct indeterminates and  $R$ be an integral domain. Set $s_1(Y)$,$\dots$,
$s_d(Y)$ for  the elementary symmetric polynomials on $Y$  and consider the polynomial
$$u=x^d-s_1(Y)x^{d-1}+\cdots +(-1)^d s_d(Y)=\prod_{1\le i \le d}(x-y_i)\ \in \ \Z[Y][x],$$ which roots are the elements in $Y$. Let $G(U)$ be the generic Gr\"obner basis of the previous remark, and define  $G(s(Y))\subset \Z[Y][X]$ as the family of polynomials obtained by specializing each variable $u_{d-k}$ into $(-1)^ks_k(Y)$ for $1\le k\le d$. \\
 Then, for any symmetric polynomial $h\in R[X]$, $r_{G(s(Y))}(h)\in R[Y][X]$ is a symmetric polynomial of degree $\le d-n$ in each variable in $X$, which is also symmetric in the elements of $Y$, and the generic symmetric Lagrange interpolant $r_Y(h)$ satisfies
 $$r_Y(h)=r_{G(s(Y))}(h)\in R[Y][X].$$
\end{corollary}
The previous corollary implies in turn  that for any multiset $A\in K$, which is a specialization of the set $Y$ of generic roots, and its associated $f\in K[x]$ and Gr\"obner basis $G$, $r_G(h)$ is in fact a polynomial in the elements of $A$:
\begin{corollary}Let $R\subset K$ be an integral domain. Given a multiset $A\subset K$ and $h\in R[X]$, we have that $r_G(h)\in R[A][ X]$.
\end{corollary}
\begin{proof}  We just evaluate the elements of $Y$ in the previous corollary into the elements of the set $A$: Then, $G=G(s(A))$  and  $r_G(h)=r_{s(Y)}(h)(A,X)\in R[A][X]$ since there are no divisions. (The order of evaluation doesn't matter since $G(s(Y))$ and $r_{s(Y)}(h)$ are symmetric in the elements of  $Y$.)
\end{proof}
In particular  this shows that in the case of a set $A$ where by Proposition~\ref{grobner} $r_A(h)=r_G(h)$,  the symmetric Lagrange interpolant  $r_A(h)$ which is defined in \ref{def:interpolant} as a rational function of the elements of the set $A$ is in fact a polynomial in them.

\smallskip
Theorem~\ref{thm:main} is then a direct consequence of the previous observations: Given a multiset $A\subset K$ and a symmetric polynomial $h\in S[X]$, the polynomial  $r_G(h)\in S_{d-n}[X]$  is well-defined and coincides with the Hermite interpolant  in the univariate case because $G=\{f\}$, and with the symmetric Lagrange interpolant   for a set $A$ by Proposition~\ref{grobner}. Therefore $r_G(h)$ is the symmetric Hermite interpolant of $h$ with respect to the multiset $A$.

%\begin{theorem}\label{thm:main3}
%Let  $K$ be a field, $A\subset K$ be a multiset as in \eqref{eq:multiset} and let $h\in S[X]$ be a symmetric polynomial.
% Then,    $$r_G(h)\in S_{d-n}[X].$$
%  Moreover,  $r_G(h)$  is a polynomial expression in the elements of $A$ which coincides with the Hermite interpolant  in the univariate case and with the symmetric Lagrange interpolant  introduced in Definition \ref{def:interpolant} in the  case of a set $A$.\\ Thus, the symmetric Hermite interpolant  $r_A(h)$ is equal to $r_G(h)$ for a multiset $A$ and $h\in S[X]$.
%\end{theorem}
%\begin{proof}  This polynomial $r_G(h)$  coincides with the Hermite interpolant  in the univariate case because $G=\{f\}$ and with the symmetric Lagrange interpolant   for a set $A$ by Proposition~\ref{grobner}.
%\end{proof}

\section{The symmetric interpolant in  Roy-Szpirglas' basis}\label{sec:RS}

The purpose of this section is to  describe the coordinates of the multivariate symmetric Hermite interpolant $r_A(h)$, defined in Theorem~\ref{thm:main},  of {\em any symmetric polynomial} $h\in S(X)$ in  the multivariate symmetric Hermite basis  presented by Roy and  Szpirglas in \cite{RoSz2020}.  For sake of completeness and to fix notations,  we first recall the construction of their basis.   Here, $\cha(K)=0$.

\subsection{Roy-Szpirglas' construction}\label{sec:RS2} {\ }

\smallskip

Note that when $A=\{a_1,\dots,a_d\}$ is a set and not a multiset, for any $A'{\subset_n}A$ we have  \begin{align*}\frac{R(X,A\setminus A')}{R(A',A\setminus A')}&=\pm \frac{v_d\big(X\cup (A\setminus A')\big)\,v_n(A')\,v_{d-n}(A\setminus A')}{v_n(X)\,v_{d-n}(A\setminus A')\,v_d(A)} \\ & =\pm \frac{v_d\big(X\cup (A\setminus A')\big)\,v_n(A')}{v_n(X)\,v_d(A)}\end{align*}
where for any (ordered) set $Y=(y_1,\dots,y_\ell)$ with $|Y|=\ell$, $$v_\ell(Y)=\prod_{1\le i<j\le \ell}(y_j-y_i)$$ is the determinant of the usual Vandermonde matrix
$$ V(Y)=\left( \begin{array}{ccc}1 & \dots & 1\\ y_1&\dots &y_\ell\\
\vdots & &\vdots\\
y_1^{\ell-1}& \dots & y_\ell^{\ell-1}\end{array}\right),$$ and $V\big(X\cup (A\setminus A')\big)$ denotes the Vandermonde matrix of the ordered set $X\cup (A\setminus A')$. 

\smallskip
In \cite{RoSz2020}, the authors  use this expression by means of Vandermonde determinants (which is related to Lagrange interpolation) to extend it to multisets, using the extension of the usual Vandermonde determinant  to  the  determinant of the  {\em generalized Vandermonde} (or {\em confluent}) matrix that arises when solving the univariate Hermite interpolation problem when there is  roots coalescence.  

\smallskip
We recall that if $A$ is a multiset as in \eqref{eq:multiset}, that we consider ordered,
the {generalized Vandermonde}
 $d\times d$ matrix
$V_d( A)$, \cite{Kal1984},
is defined as
$$
V_d(
A)=\begin{array}{|c|c|c|c}
\multicolumn{3}{c}{\scriptstyle{d}}\\
\cline{1-3} & & & \\ V_d(a_1,d_1) & \dots &V_d(a_m,d_m)
 &\scriptstyle d\\ & & & \\
 \cline{1-3}
 \multicolumn{2}{c}{}
\end{array},
$$
where
$$
V_d(a_i,d_i):=\begin{array}{|cccccc|c}
\multicolumn{6}{c}{\scriptstyle{d_i}}\\
\cline{1-6}
1&0 &0&0&\dots & 0 &\\
a_i&1 &0&0&\dots & 0 &\\
a^2_i&{2\choose 1}a_i &1&0&\dots & 0 &\\[.5mm]
a^3_i&{3\choose 1}a_i^2 &{3\choose 2}a_i&1&\dots & 0 &\\
\vdots&  \vdots& \vdots & &\ddots &\vdots &\scriptstyle d\\
\vdots&  \vdots& \vdots & &  &1 &\\
\vdots&  \vdots& \vdots & &   &\vdots  &\\
a_i^{d-1}&{d-1\choose 1}a_i^{d-2} &{d-1\choose 2}a_i^{d-3}&\dots &\dots & {d-1\choose d_i-1}a_i^{d-d_i}& \\
 \cline{1-6}
 \multicolumn{2}{c}{}
\end{array}.
$$
Its determinant $v_d(A)$ satisfies, \cite{Ait1939},
\begin{equation*}\label{gVd}
v_d(A)=\prod_{1\leq i<j\leq m}
(a_j-a_i)^{d_id_j}.
\end{equation*}
Note that matrix $V_d(A)$ is the matrix which arises when solving the Hermite linear system
$$\left\{\begin{array}{l}r(a_1)={c_{1,0}},\ r'(a_1)=1!\,{c_{1,1}},\ \dots\ ,\ r^{(d_1-1)}(a_1)={(d_1-1)!}\,{c_{1,d_1-1}}  \\
\ \ \vdots \\
r(a_m)={c_{m,0}},\ r'(a_m)=1!\,{c_{m,1}},\ \dots \ ,\ r^{(d_m-1)}(a_m)={(d_m-1)!}\,{c_{m,d_m-1}}, \end{array}\right. ,
$$
where the multiplicative factors $k!$ are added to avoid a multiplicative factor  in $v_d(A)$.
For example, for $A:=\{a,a,a;b,b\}$,
\begin{equation}\label{eq:v5}V_5(A)=\left[\begin{array}{ccc|cc}
1&0&0&1&0\\
a&1&0&b&1\\
a^2&2a&1&b^2&2b\\
a^3&3a^2&3a&b^3&3b^2\\
a^4&4a^3& 6a^2&b^4& 4b^3
                       \end{array}
\right] \quad \mbox{with} \quad  v_5(A)=(b-a)^6,\end{equation}
is the invertible matrix which arises when looking for a polynomial of degree at most 4 which satisfies the 5 conditions
$$r(a)=c_0,\ r'(a)=c_1,\  r''(a)=2\,c_2, \ r(b)=d_0,\ r'(b)=d_1 ,$$
and by Cramer's rule, one can see that the Hermite basis $\mathcal B$ can be described by means of signed minors of $V_5(A)$ divided by $v_5(A)$. 

\smallskip
When dealing with a multiset $A$ as in \eqref{eq:multiset},  that we consider ordered, and its corresponding generalized Vandermonde matrix  $V_d(A)$, it is convenient to label the  elements of $A$ in the following way:
\begin{equation}\label{eq:multiset2}A=\big(a_{1,0},\dots,a_{1,d_1-1};\dots;a_{m,0},\dots,a_{m,d_m-1}\big)\end{equation}
 {where} $a_{i,j}:=a_i$ for $1\le i\le m$, $0\le j< d_i$, in order to associate to each element in $A$ a column in $V_d(A)$. \\ In \cite{RoSz2020}, the authors interpret     each (ordered) ``subset" $A'{\subset_n}A$ as a choice of columns in  $V_d(A)$ indexed by the elements of $A'$, and    $V(A\setminus A')\in K^{d\times (d-n)}$ corresponds to the submatrix of $V_d(A)$ obtained by deleting the columns indexed by $A'$. \\
 For instance, for our  previous example  $A=\{a,a,a;b,b\}=\big( a_0,a_1,a_2;b_0,b_1\big)$, with $a_i=a$ and $b_j=b$, for $A_1=\big(a_1;b_0\big)\subset_{2}  A$, $V(A\setminus A_1)$ is obtained by deleting the second column of $V_5(A)$ in \eqref{eq:v5}, corresponding to the second element of the $a$-block, and its fourth column,  corresponding to the first element of the $b$-block: We obtain
 $$V(A\setminus A_1)=\left[\begin{array}{cc|c}
1&0&0\\
a&0&1\\
a^2&1&2b\\
 a^3&3a&3b^2\\
 a^4& 6a^2& 4b^3
                       \end{array}
\right].$$
 Finally, given  $X=(x_1,\dots,x_n)$, we define $V(X\cup (A\setminus A'))$ as  the matrix obtained by attaching to the left of $V(A\setminus A')$ the usual   Vandermonde matrix  of size $d\times n$ of the $n$ ordered variables in $X$. In our example where $X=(x,y)$ this gives the $5\times 5$ matrix
$$V\big(X\cup(A\setminus A_1)\big)=\left[\begin{array}{cc|cc|c}
1&1&1&0&0\\
x&y&a&0&1\\
x^2&y^2&a^2&1&2b\\
x^3&y^3& a^3&3a&3b^2\\
x^4&y^4& a^4& 6a^2& 4b^3
                       \end{array}
\right].$$

The following result is
 \cite[Prop.18]{RoSz2020}:
 \begin{theorem}\label{thm:RS} {\em (Symmetric Hermite basis.)} \\  The set
$$\mathcal{B}:=\Big\{\  \frac{\det\big(V\big(X\cup (A\setminus A')\big)\big)}{v_n(X)}:\ A'{\subset_n}A\ \Big\}$$ is a basis for $S_{d-n}(X)$, the set of all symmetric polynomials in $n$ variables of degree in each variable bounded by $d-n$, which  extends the symmetric Lagrange basis $\{\,{R(X,A\setminus A')}:\ A'{\subset_n}A \,\} $ to the case of multisets, as well as the usual basis for the univariate Hermite interpolation.
\end{theorem}
To describe the coordinates of a symmetric polynomial $r\in S_{d-n}[X]$ in this basis as they do, we need to introduce the following evaluation and derivative of a polynomial with respect to $A'$:\\
 Given any polynomial $g\in K[X]$ (nonnecessarily symmetric)  in the ordered variables $X=(x_1,\dots,x_n)$ and given $ A'{\subset_n}A$
 as in \eqref{eq:multiset2},
\begin{itemize}
\item
$g(A')\in K$ means evaluating for all $1\le k\le n$ the variable $x_k$ in $a_i$ if the $k-$th element of $A'$ is $a_{i,j}$ for some $0\le j<d_i$.
\item $\partial^{A'}g\in K[X]$ means that for each $1\le k\le n$, if the $k-$element of $A'$ is $a_{i,j}$ for $1\le i\le m$ and $0\le j< d_i$, we derivate
 $j$ times the variable $x_k$ in $g$ and divide by $j!$.
\end{itemize}
To continue with our example where $A=\{a,a,a;b,b\}$,  for $A_1=(a_1;b_0), A_2=(a_0;b_0)\,\subset_{2}  A$,  and $$v_1:=\det\big(V\big(X\cup (A\setminus A_1)\big)\big),\ v_2:=\det\big(V\big(X\cup (A\setminus A_2)\big)\big),$$
we get $v_1(A_1)=0$ since there are two repeated columns in the matrix $V\big(X\cup (A\setminus A_1)\big)$, while $v_2(A_2)=(b-a)^6$ since it is the determinant of a matrix that by 2 transpositions of columns gives $V_5\big(A)$.  Also, $\partial^{A_1}(v_1)(A_1)= -(b-a)^6$
since it is the determinant of the following matrix
$$\left[\begin{array}{cc|cc|c}
0&1&1&0&0\\
1&y&a&0&1\\
2x&y^2&a^2&1&2b\\
3x^2&y^3& a^3&3a&3b^2\\
4x^3&y^4& a^4& 6a^2& 4b^3
                       \end{array}
\right]$$ evaluated into $x=a$ and $y=b$, which by 3 transpositions of columns gives the matrix $V_5(A)$ and $\partial^{A_1}(v_2)(A_1)= 0$ since again there are two repeated columns in the evaluation of the matrix to be considered.

\smallskip
In  \cite[Prop.19]{RoSz2020},
the coordinates $c_{A'}$ in the basis  $\mathcal{B}$ of any polynomial $r\in S_{d-n}[X]$, that is of any polynomial that is already of degree $\le d-n$ in each variable, are described as follows:
\begin{proposition}\label{eq:coeff} {\em (Coordinates in the Hermite basis.)}
$$r(X)=\sum_{A'{\subset_n}A} c_{A'} \frac{\det\big(V\big(X\cup (A\setminus A')\big)\big)}{v_n(X)}$$
with ~\begin{equation*}c_{A'}= \varepsilon_{A'}\frac{\partial^{A'}(v_n\, r)(A')}{v_d(A)}\ \in K,\end{equation*}
where $\varepsilon_{A'}$ is $\pm 1$  according to  $\partial^{A'}\det\big(V\big(X\cup (A\setminus A')\big)\big)(A')=\pm v_d(A)$.
\end{proposition}

\noindent  {\bf Example.} Let us illustrate Theorem~\ref{thm:RS} and Proposition~\ref{eq:coeff}  with  a smaller example: $$S_2[x,y] \ \mbox{ and } \ A=(a,a;b,b)=(a_0,a_1;b_0,b_1) \ \mbox{ with } \ a_1=a_2=a \ , \ b_1=b_2=b.$$
Here $v_2(A)=y-x$.
We need to consider the following 6 ``subsets" of $A$,
\begin{align*}A_1=(a_0,a_1),& \ A_2= (a_0;b_0), \ A_3=(a_0;b_1),\\ &
A_4= (a_1;b_0),\ A_5=(a_1;b_1),\ A_6= (b_0,b_1),\end{align*}
with  their corresponding matrices $V_i:=V(\{x,y\}\cup(A\setminus A_i))$ and  determinants $v_i:=\det(V_i)$. We set $\omega_i=v_i/(y-x)$, $1\le i\le 6$.
%\begin{align*}&A_1=(a_0,a_1), \quad V_1:=V(\{x,y\}\cup(A\setminus A_1))=\left[\begin{array}{cc|cc}
%1& 1& 1&0\\
%x&y& b&1\\
%x^2&y^2&b^2&2b\\
%x^3&y^3& b^3& 3b^2                    \end{array}
%\right],\\
%&A_2= (a_0;b_0),\quad V_2:=V(\{x,y\}\cup(A\setminus A_2))=\left[\begin{array}{cc|c|c}
%1& 1& 0&0\\
%x&y& 1&1\\
%x^2&y^2&2a&2b\\
%x^3&y^3&   3a^2& 3b^2                  \end{array}
%\right],\\ &
%A_3=(a_0;b_1), \quad V_3:=V(\{x,y\}\cup(A\setminus A_3))=\left[\begin{array}{cc|c|c}
%1& 1& 0&1\\
%x&y& 1&b\\
%x^2&y^2&2a&b^2\\
%x^3&y^3&   3a^2& b^3                  \end{array}
%\right] ,\\ &
%A_4= (a_1;b_0), \quad V_4:=V(\{x,y\}\cup(A\setminus A_4))=\left[\begin{array}{cc|c|c}
%1& 1& 1&0\\
%x&y& a&1\\
%x^2&y^2&a^2&2b\\
%x^3&y^3&   a^3& 3b^2                  \end{array}
%\right] ,\\ &
%A_5=(a_1;b_1), \quad V_5:=V(\{x,y\}\cup(A\setminus A_5))=\left[\begin{array}{cc|c|c}
%1& 1& 1&1\\
%x&y& a&b\\
%x^2&y^2&a^2&b^2\\
%x^3&y^3&   a^3& b^3                  \end{array}
%\right],\\
%&A_6= (b_0,b_1), \quad V_6:=V(\{x,y\}\cup(A\setminus A_6))=\left[\begin{array}{cc|c|c}
%1& 1& 1&0\\
%x&y& a&1\\
%x^2&y^2&a^2&2a\\
%x^3&y^3&a^3&   3a^2                  \end{array}
%\right].
%\end{align*}
%Setting $v_i:=\det(V_i)$, we easily see that $\dfrac{v_i}{y-x}\in S_2[x,y]$.\\
Then the  set
$$\mathcal{B}=\big\{\omega_1, \dots,\omega_6 \big\}$$ is the symmetric Hermite interpolation basis for $S_2[x,y]$,  and
any symmetric polynomial $r\in S_2[x,y]$ can be written  as
$$r(x,y)= \sum_{i=1}^6 c_i \omega_i$$
where  for $1\le i\le 6$, $$c_i=\varepsilon_i\, \frac{\partial^{A_i} \big( (y-x)r\big)(A_i)}{(b-a)^4}  \quad \mbox{with} \quad \varepsilon_i=-1  \Leftrightarrow i=2 \ \mbox{or} \ 5.$$
%This happens because
%\begin{align*}
%&\partial^{A_1}v_1(A_1)= \partial_yv_1(a,a)=(b-a)^4 \ \mbox{ and } \ \partial^{A_1}(v_j)(A_1)= 0 \mbox{ for } j\ne 1,\\
%&\partial^{A_2}v_2(A_2)= v_2(a,b)=-(b-a)^4 \ \mbox{ and } \ \partial^{A_2}(v_j)(A_2)= 0 \mbox{ for } j\ne 2,\\
%&\partial^{A_3}v_3(A_3)= \partial_yv_3(a,b)=(b-a)^4 \ \mbox{ and } \ \partial^{A_3}(v_j)(A_3)= 0 \mbox{ for } j\ne 3,\\
%&\partial^{A_4}v_4(A_4)= \partial_xv_4(a,b)=(b-a)^4 \ \mbox{ and } \ \partial^{A_4}(v_j)(A_4)= 0 \mbox{ for } j\ne 4,\\
%&\partial^{A_5}v_5(A_5)= \partial_{xy}v_5(a,b)=-(b-a)^4 \ \mbox{ and } \ \partial^{A_5}(v_j)(A_5)= 0 \mbox{ for } j\ne 5\\
%&\partial^{A_6}v_6(A_6)= \partial_yv_6(b,b)=(b-a)^4 \ \mbox{ and } \ \partial^{A_6}(v_j)(A_16)= 0 \mbox{ for } j\ne 6.
%\end{align*}
For information we have:
\begin{align*} &\omega_1=  (x-b)^2(y-b)^2 \quad \mbox{and} \quad \omega_6=(x-a)^2(y-a)^2 , \\
& \omega_2= (b-a)\big(2(x^2+xy+y^2) - 3(a+b)(x+y) + 6ab\big),\\
&\omega_3=-(x-b)(y-b)\big(xy+(-2a+b)(x+y)+3a^2-2ab\big), \\&
\omega_4=(x-a)(y-a)\big(xy+(a-2b)(x+y)+3b^2-2ab\big),\\
& \omega_5= (b-a)(x-a)(x-b)(y-a)(y-b).
\end{align*}
Setting $g=(y-x)r$ for $r\in S_2[x,y]$ we then obtain
\begin{align*} r & =\frac{1}{(b-a)^4} \Big( \partial_y g(a,a) \omega_1 - g(a,b)\omega_2 + \partial_y g(a,b) \omega_3  \\ & \qquad \qquad \qquad \quad +\partial_x g(a,b) \omega_4- \partial_{xy} g(a,b) \omega_5 + \partial_y g(b,b) \omega_6\Big). \end{align*}

\subsection{The coordinates of any symmetric polynomial in $\mathcal{B}$}\label{sec:coord} {\ }

\smallskip

In this section we extend Proposition~\ref{eq:coeff} to any symmetric polynomial $h\in S[X]$, of any degree: We show that the description of the coordinates of the polynomial $r_h \in S_{d-n}[X]$ in the basis $\mathcal{B}$ coincides in terms of $h\in S[X]$ with the description of the coordinates of $r\in S_{d-n}[X]$ given in that proposition. To this aim, we present the following
result that uses a connection between the normal form $r_G(h)$ of a symmetric polynomial $h\in S[X]$  and its normal form $r_F(h)$ modulo the Gr\"obner basis  $F:=\{f(x_1), \ldots, f(x_{n })\}$ of the ideal
$\big( f(x_1),\dots,f(x_n)\big) \subset K[X]$  w.r.t. the lexicographic order $\prec$.  The idea of the proof  is similar to \cite[Lemma 11]{KoRo2012}.

\begin{proposition}\label{prop:Vand}
Let $A$ be a multiset as in \eqref{eq:multiset} and
$f= \prod_{1\le i\le m}(x-a_i)^{d_i} $.
 Denote $$F=\{ f(x_1), \ldots, f(x_{n})\} \quad \mbox{and}\quad
v_{n}:=\prod_{1\leq i<j\leq n} (x_j-x_i)
$$
for the Vandermonde polynomial in the $n$ (ordered) variables in $X$. Then, for any symmetric polynomial  $h\in S[X]$ we have
$$r_A(h)= \frac{r_F(v_{n}\, h)}{v_{n}}.
$$
\end{proposition}

\begin{proof} In view of Theorem~\ref{thm:main},
since $r_A(h)=r_G(h)$, where $G$ is the symmetric Gr\"obner basis introduced in \eqref{GB}, it is enough to show that $$ {r_F(v_{n}\, h)}={v_{n}}\, r_G(h).$$
We first prove  Identity \eqref{eq:fg} by recurrence on  $1\le k\le n$:\\
The case $k=1$ is by definition.\\
Now, let $1\le k<n$. By the inductive hypothesis, the $k$-th identity in \eqref{eq:fg} applied to $x_{k+1}$ instead of $x_k$ reads
\begin{align*}f(x_{k+1})&=\prod_{i=1}^{k-1}(x_{k+1}-x_i)\, g_k(x_1,\dots,x_{k-1},x_{k+1}) + \prod_{i=1}^{k-2}(x_{k+1}-x_i)\, g_{k-1}(x_1,\dots,x_{k-1}) \\ & \qquad \qquad +\cdots + (x_{k+1}-x_1)\, g_2(x_1,x_2) + g_1(x_1).\end{align*}
The claim for $k+1$ thus follows by applying
\eqref{eq:idg}, which  says that the symmetric polynomials $g_k$ and $g_{k+1}$ satisfy
$$g_{k}(x_1,\dots,x_{k-1},x_{k+1})= (x_{k+1}-x_k)\,g_{k+1}(x_1,\dots,x_{k+1})+g_k(x_1,\dots,x_k).$$
Identity \eqref{eq:fg} implies the following inclusions of ideals
$$
(F) \subset (G)\subset \big(( F ) :   v_{n}\big),
$$
where $\big(( F ) :   v_{n}\big)=\{g\in K[X]: v_n\,g\in (F) \} $ is the colon ideal, and in particular $v_n\,g_i\in (F)$ for $1\le i\le n$.
 Let
$$
h=\sum_{i=1}^{n} p_i\,g_{i} +r_G(h)\quad \in \ S[X]
$$ for some $p_1, \ldots, p_{n}\in K[X]$ and $r_G(h)\in S_{d-n}[X]$,  since $h$ is symmetric.
Therefore
$$
v_{n}\, h =  \sum_{i=1}^{n} p_i \,v_{n} \,g_{i}  + v_{n}\, r_G(h) \ = \
\sum_{i=1}^{n} q_i\,f(x_i) + v_{n}\, r_G(h),
$$
for some $q_1, \ldots, q_{n}\in K[X]$. \\We note that
$v_{n}\,r_G(h)$ is reduced modulo the Gr\"obner basis $F$, since the degree
of $v_{n}$ in each variable is $n-1$  and the max-degree of $r_G(h)$ is at most $d-n$, so $v_{n}\, r_G(h)$ has degree   at most $d-1$ in each variable.
Thus $$v_{n}\, r_G(h)=r_F(h),$$ which implies
the statement.
\end{proof}

\begin{corollary}\label{cor:Hermiteinterpolant}{\em (Coordinates of the symmetric Hermite interpolant.)}\\
Given a multiset $A$ as in \eqref{eq:multiset}, with $|A|=d\ge n$, and given any symmetric polynomial $h\in S[X]$, the {\em symmetric Hermite interpolant}  $r_A(h)\in S_{d-n}[X]$ of $h$ with respect to the multiset $A$ satisfies
$$r_A(h)=\sum_{A'{\subset_n}A} c_{A'} \frac{\det\big(V\big(X\cup (A\setminus A')\big)\big)}{v_n(X)}$$
with ~\begin{equation*}c_{A'}= \varepsilon_{A'}\frac{\partial^{A'}(v_n\, h)(A')}{v_d(A)}\ \in K,\end{equation*}
where $\varepsilon_{A'}$ is $\pm 1$  according to  $\partial^{A'}\det\big(V\big(X\cup (A\setminus A')\big)\big)(A')=\pm v_d(A)$.
%$$r_A(h)=\frac{1}{v_d(A)v_n}\sum_{A'{\subset_n}A}\varepsilon_{A'}\partial^{A'}(v_n\, h)(A') \det\big(V\big(X\cup (A\setminus A')\big)\big),$$
%where $\varepsilon_{A'}$ is $\pm 1$  according to  $\partial^{A'}\det\big(V\big(X\cup (A\setminus A')\big)\big)(A')=\pm v_d(A)$.
\end{corollary}

\begin{proof}
By Proposition~\ref{prop:Vand},
$$v_n\, h=\sum_{i=1}^{n}  q_i\, f(x_i) + v_n\,r_A(h) $$
for some $q_i\in K[X]$ and therefore
$$\partial^{A'}(v_n\, h) (A')=  \partial^{A'}\big(v_n\, r_A(h)\big)(A') \quad \mbox{for all}\quad A'{\subset_n}A.$$
We conclude by applying Proposition~\ref{eq:coeff}.
\end{proof}

\section{The Taylor interpolant} \label{sec:taylor}
As an application we  study the case when $$ A=\{\underbrace{0,\dots 0}_{d}\}= \big(0_0,\dots, 0_{d-1}\big)$$
which corresponds with the interpolation polynomial given by the Taylor expansion at $a=0$, and relate it to the  Schur polynomials.

\smallskip
In this case, we need to consider for $A'=\big(0_{i_1},\dots 0_{i_n}\big){\subset_n}A$ with $0\le i_1< \cdots <i_n<d$ the determinant of the submatrix of
$$V\big(X\cup A\big) =
\left(\begin{array}{ccc|cccc}1&\dots&1&1&  & & \\
x_1&\dots &x_n& &1  & & \\
\vdots & & \vdots && & \ddots & \\
x_1^{d-1}& \dots & x_n^{d-1}& & & & 1
\end{array}
\right) \quad \in \ K^{d\times (n+d)}
 $$
 where the columns indexed by the columns corresponding to $A'$ have been deleted.
 Therefore,
 for each of these $A'$, $$\frac{\det\big(V\cup (A\setminus A')\big)}{v_n(X)}=\pm s_{A'}(X)$$ where $ s_{A'}(X)$ denotes the Schur polynomial indexed by the rows corresponding to $A'$.

\smallskip
The previous section shows (in an alternative way) that the set of Schur polynomials
 $$\{ s_{A'}\,:\, A'{\subset_n}  A\}$$
 is a basis of $S_{d-n}(X)$, the vector space of all symmetric polynomials in $n$ variables of degree in each variable $\le d-n$: this is a known  result that can be found for instance in  \cite[Claim (3.2)]{Mac1995}. Moreover, the Hermite interpolant of a symmetric polynomial $h\in S[X]$ with respect to this multiset $A$ is given by
 $$r_A(h)=\sum_{A'{\subset_n}A} \partial^{A'}(v_n\,h)(A') s_{A'},$$
 because the sign $\varepsilon_{A'}$ is canceled by the sign of $\det\big(X\cup(A\setminus A')\big)/v_n(X)$ with respect to $s_n(X)$.

\medskip
\noindent {\bf Example.}
For $d=4$, $A=\{0,0,0,0\}=(0_0,0_1,0_2,0_3)$,  and $n=2$. The Schur polynomials are $1, x+y, x^2+xy+y^2, xy, x^2y+xy^2$ and $x^2y^2$.\\  For $h\in S[x,y]$ and $g:=(y-x)h$, setting ${\bfzero}=(0,0,0,0)$, we have that
\begin{align*}
r_A(h)& = \frac{\partial g}{\partial y}(\bfzero) + \frac{\partial^2 g}{2\partial y^2}(\bfzero)(x+y) +\frac{\partial^3 g}{6\partial y^3}(\bfzero) (x^2+xy+x^2)\\ & \quad  + \frac{\partial^3 g}{2\partial x \partial y^2}(\bfzero)  xy + \frac{\partial^4 g}{6\partial x \partial y^3}(\bfzero)  (x^2y+xy^2) + \frac{\partial^5 g}{12\partial x^2\partial y^3}(\bfzero)  x^2y^2\\
&= h(\bfzero) + \frac{\partial h}{\partial y}(\bfzero)(x+y)+ \frac{\partial^2 h}{2\partial y^2}(\bfzero)(x^2+xy+x^2)\\
&\quad  + \big(\frac{\partial^2 h}{\partial x\partial y}(\bfzero)-\frac{\partial^2 h}{2\partial y^2}(\bfzero)\big)xy  +
\big(\frac{\partial^3 h}{2\partial x\partial y^2}(\bfzero)-\frac{\partial^3 h}{6\partial y^3}(\bfzero)\big) (x^2y+xy^2)\\ &
\quad +  \big(\frac{\partial^4 h}{4\partial x^2\partial y^2}(\bfzero)-\frac{\partial^4 h}{6\partial x\partial y^3}(\bfzero)\big) x^2y^2.
\end{align*}

\bigskip
\bibliographystyle{alpha}
\def\cprime{$'$} \def\cprime{$'$} \def\cprime{$'$}

\end{document}